\documentclass[11pt,twoside]{article}

\usepackage{amsmath,amssymb,amsthm,url}

\newcommand{\Prim}{\mathrm{Prim}}

\newcommand{\vol}{\mathrm{vol}}

\newcommand{\as}{\quad\text{as}\quad}
\newcommand{\tinf}{\to\infty}
\newcommand{\disp}{\displaystyle}
\newcommand{\bsla}{\backslash}

\newcommand{\cD}{\mathcal{D}}

\newcommand{\bC}{\mathbb{C}}
\newcommand{\bR}{\mathbb{R}}

\newcommand{\bZ}{\mathbb{Z}}

\newcommand{\noi}{\noindent}
\renewcommand{\Re}{\mathrm{Re}}
\renewcommand{\Im}{\mathrm{Im}}

\newcommand{\divset}{\hspace{3pt}|\hspace{3pt}}

\newcommand{\gam}{\gamma}
\newcommand{\Gam}{\Gamma}

\newcommand{\sr}{\mathrm{SL}_2(\bR)}
\newcommand{\sz}{\mathrm{SL}_2(\bZ)}

\newcommand{\vcpt}{\vspace{12pt}}

\newtheorem{thm}{Theorem}[section]
\newtheorem{prop}[thm]{Proposition}
\newtheorem{lem}[thm]{Lemma}

\numberwithin{equation}{section}

\title{Selberg's zeta function for the modular group in the critical strip}
\author{Yasufumi Hashimoto}
\date{}

\begin{document}
\markboth
{Y. Hashimoto}
{Selberg's zeta functions for modular group}
\pagestyle{myheadings}

\maketitle
\renewcommand{\thefootnote}{}

\begin{abstract}
In the present paper, we study the growth of the Selberg zeta function for the modular group
in the critical strip.
\end{abstract}

\section{Introduction and the main theorem}
Let $H:=\{x+y\sqrt{-1}\divset x,y\in\bR,y>0\}$ be the upper half plane 
and $\Gam$ a discrete subgroup of $\sr$ with $\vol(\Gam\bsla H)<\infty$. 
Denote by $\Prim(\Gam)$ the set of primitive hyperbolic conjugacy classes 
of $\Gam$ and $N(\gam)$ the square of the larger eigenvalue of $\gam$. 
The Selberg zeta function for $\Gam$ is defined by 
\begin{align*}
Z_{\Gam}(s):=\prod_{\gam\in\Prim(\Gam),n\geq0}
(1-N(\gam)^{-s-n}),\qquad \Re{s}>1.
\end{align*}
It is well known that $Z_{\Gam}(s)$ is analytically continued to the whole complex plane 
and has a functional equation between the values at $s$ and $1-s$ (see, e.g. \cite{He}). 
The aim of the present paper is to study the growth of $Z_{\Gam}(s)$ 
in the critical strip as $|\Im{s}|\tinf$.

For the Riemann zeta function 
$$\zeta(s):=\prod_{p}(1-p^{-s})^{-1}, \qquad \Re{s}>1,$$
there have been various works on the growth in the critical strip. 
In fact, it was proved \cite{Hux} that  
$$\zeta(1/2+iT)\ll_{\epsilon}T^{\frac{32}{205}+\epsilon}$$ 
and has been considered that 
$\zeta(1/2+iT)\ll_{\epsilon}T^{\epsilon}$ as $T\tinf$. 
On the other hand, 
it is known that 
\begin{align}
\log{Z_{\Gam}(\sigma+iT)}, \frac{Z'_{\Gam}(\sigma+iT)}{Z_{\Gam}(\sigma+iT)}
\ll_{\epsilon}T^{2-2\sigma+\epsilon}, 
\as T\tinf \label{SZorg}
\end{align}
for $1/2<\sigma<1$ (see \cite{He} for co-compact $\Gam$ and \cite{IwPGT} for $\Gam=\sz$). 
This means that the growth of $Z_{\Gam}(s)$ is exponential of $|\Im{s}|$ 
and is quite different to the Riemann zeta function. 

In the present paper, we improve \eqref{SZorg} for the modular group.
The main result is as follows.

\begin{thm} \label{thm} 
Let $\Gam=\sz$ and $s=\sigma+iT$ with $1/2<\sigma<1$. 
Then we have  
\begin{align}\label{thmeq}
\frac{Z'_{\Gam}(s)}{Z_{\Gam}(s)}\ll_{\epsilon}
\begin{cases}
\disp T^{\frac{19}{9}-\frac{20}{9}\sigma+\epsilon}, & (\frac{1}{2}<\sigma\leq \frac{5}{8}),\\
\disp T^{\frac{52}{27}(1-\sigma)+\epsilon}, & (\frac{5}{8}< \sigma<1),
\end{cases} \as T\tinf
\end{align}
for any $\epsilon>0$.
\end{thm} 

To prove Theorem \ref{thm}, we first describe $Z'_{\Gam}(s)/Z_{\Gam}(s)$ 
by a sum over $\gam\in\Prim(\Gam)$ as given in Proposition \ref{prop1}. 
For $\Gam=\sz$, such a sum can be expressed in terms of 
the class numbers and the fundamental units 
of the primitive indefinite binary quadratic forms, 
and then 
we can prove Theorem \ref{thm} 
by using the large sieve method for the Dirichlet $L$ function 
and van-der Corput's exponential sum estimate.

\section{Proof of Theorem \ref{thm}}

We first state the following explicit formula 
for $Z'_{\Gam}(s)/Z_{\Gam}(s)$ to prove Theorem \ref{thm}.
\begin{prop}\label{prop1}
Let $s=\sigma+iT\in \bC$ with $1/2<\sigma<1$.  
Then we have 
\begin{align*}
\frac{Z'_{\Gam}(s)}{Z_{\Gam}(s)}=
&\sum_{\begin{subarray}{c}\gam\in\Prim(\Gam),j\geq1 \\ N(\gam)^j<x\end{subarray}} 
\frac{\log{N(\gam)}}{1-N(\gam)^{-j}}
\left(1-\frac{N(\gam)^{j}}{x} \right)N(\gam)^{-js} \\
&+O_{\epsilon}\left(T^{-2}x^{1-\sigma}+T^{1+\epsilon}x^{\frac{1}{2}-\sigma+\epsilon}\right), 
\as T,x\tinf
\end{align*}
for any $\epsilon>0$.
\end{prop}

\noi{\it Proof.} 
For $V=T^{3}$ and $\epsilon>0$, let $C$ be the rectangle 
with the corners $1+\epsilon-iV$, $1+\epsilon+iV$, $1/2+\epsilon+iV$, $1/2+\epsilon-iV$, 
and 
\begin{align*}
J:=\frac{1}{2\pi i}\int_{\partial C}\frac{Z'_{\Gam}(z)}{Z_{\Gam}(z)}\frac{x^{z-s}}{(z-s)(z-s+1)}dz,
\end{align*}
where the integral is in an anti-clockwise direction.
Since the singular point of $Z'_{\Gam}(z)/Z_{\Gam}(z)$ in $C$ is only 
a single pole at $z=1$, 
we have 
\begin{align}
J=\frac{Z'_{\Gam}(s)}{Z_{\Gam}(s)}+\frac{x^{1-s}}{(1-s)(2-s)}
\end{align}
by the residue theorem.
Next, divide $J$ by 
\begin{align*}
2\pi i J=\int_{\partial C}=&\int_{1+\epsilon-i\infty}^{1+\epsilon+i\infty}
-\int_{1+\epsilon-i\infty}^{1+\epsilon-iV}-\int_{1+\epsilon+iV}^{1+\epsilon+i\infty}
+\int_{1+\epsilon+iV}^{1/2+\epsilon+iV}+\int_{1/2+\epsilon+iV}^{1/2+\epsilon-iV}
+\int_{1/2+\epsilon-iV}^{1+\epsilon-iV}\\
=:&J_{1}-J_{2}-J_{3}+J_{4}+J_{5}+J_{6}.
\end{align*}
Due to \eqref{SZorg}, we can bound $J_{2},\dots,J_{6}$ by 
\begin{align}
J_{2},J_{3}
&\ll_{\epsilon}\int_{V}^{\infty}|u|^{\epsilon}\frac{x^{1-\sigma}}{(u-T)^{2}+1}du
\ll V^{-1+\epsilon}x^{1-\sigma},\\
J_{4},J_{6}
&\ll_{\epsilon} V^{1+\epsilon}x^{1-\sigma}(V-T)^{-2}\ll V^{-1}x^{1-\sigma},\\
J_{5}&\ll_{\epsilon} \int_{-V}^{V}|u|^{1+\epsilon}\frac{x^{1/2-\sigma}}{(u-T)^{2}+1}du
\ll_{\epsilon} T^{1+\epsilon}x^{1/2-\sigma}.
\end{align}
The remaining term $J_{1}$ is given by 
\begin{align}
\frac{1}{2\pi i}J_{1}=& \sum_{\begin{subarray}{c}\gam\in\Prim(\Gam),j\geq1 \end{subarray}} 
\frac{\log{N(\gam)}}{1-N(\gam)^{-j}}N(\gam)^{-js}
\cdot \frac{1}{2\pi i}\int_{1-s+\epsilon-i\infty}^{1-s+\epsilon+i\infty} 
\frac{(x/N(\gam)^{j})^{z}}{z(z+1)}dz\notag\\
=&\sum_{\begin{subarray}{c}\gam\in\Prim(\Gam),j\geq1 \\ N(\gam)^j<x\end{subarray}} 
\frac{\log{N(\gam)}}{1-N(\gam)^{-j}}
\left(1-\frac{N(\gam)^{j}}{x} \right)N(\gam)^{-js}.
\end{align} 
The proposition follows from (2.1)--(2.5). \qed

\vcpt

It is known that there is a one-to-one correspondence 
between equivalence classes of primitive indefinite binary quadratic forms 
and the elements of $\Prim(\sz)$, 
and its correspondence is given as follows.
\begin{align}\label{1to1}
[a,b,c]=ax^2+bxy+cy^2
\leftrightarrow
\begin{pmatrix}\disp\frac{t+bu}{2}&-cu\\ 
au& \disp\frac{t-bu}{2}\end{pmatrix},
\end{align}
where $D:=b^{2}-4ac$ is the discriminant of $[a,b,c]$ 
and $(t,u)$ is the smallest positive solution of 
the Pell equation $t^{2}-Du^{2}=4$ (see e.g. Chap. 5 in \cite{Ga}).
Due to the correspondence above, 
we have  
\begin{align*} 
\psi_{s}(x):=&\sum_{\begin{subarray}{c}\gam\in\Prim(\Gam),j\geq1, \\ N(\gam)^j<x \end{subarray}} 
\frac{\log{N(\gam)}}{1-N(\gam)^{-j}}
\left(1-\frac{N(\gam)^{j}}{x} \right)N(\gam)^{-js}\\
=&\sum_{\begin{subarray}{c}D\in\cD,j\geq1\\ \epsilon_{1}(D)^{j}<x^{1/2}\end{subarray}}
\frac{2\log{\epsilon_{1}(D)}h(D)}{1-\epsilon_{1}(D)^{-2j}}
\left(1-\frac{\epsilon_{1}(D)^{2j}}{x}\right)\epsilon_{1}(D)^{-2js},
\end{align*}
where $\cD$ is the set of $D>0$ with $D\equiv 0,1\bmod{4}$, 
$\epsilon_{1}(D):=\frac{1}{2}(t+u\sqrt{D})=\frac{1}{2}(t+\sqrt{t^{2}-4})$
and $h(D)$ is the class number in the narrow sense 
(see also \cite{Sa1}).
Applying the works of Kuznetzov and Bykovskii \cite{Kuz,Byk} (see also \cite{SY,BF}), 
we can write $\psi_{s}(x)$ as follows.
\begin{align}\label{psi} 
\psi_{s}(x)=&2\sum_{3\leq t<X}L(1,t^{2}-4)\left(1-\frac{\epsilon(t)^{2}}{x}\right)
\epsilon(t)^{1-2s},
\end{align}
where $X:=x^{1/2}+x^{-1/2}$, 
$\epsilon(t):=\frac{1}{2}(t+\sqrt{t^{2}-1})$ and
\begin{align*}
L(z,D):=\sum_{dl^{2}=D}l^{1-2z}L\left(z,\left(\frac{d}{*}\right)\right)=\sum_{dl^{2}=D}l^{1-2z}
\sum_{n\geq1}\left(\frac{d}{n}\right)n^{-z}.
\end{align*}

To estimate $\psi_{s}(X)$, we prepare the following two lemmas.

\begin{lem}\label{lemHB} (Theorem 2 of \cite{HB})
For $z\in\bC$ and a Dirichlet character $\chi$, let 
$L(z,\chi):=\sum_{n\geq1}\chi(n)n^{-z}$ 
be the Dirichlet $L$-function. 
Then we have 
\begin{align*}
\sum_{\chi\in S(Q)}\left| L\left(\frac{1}{2}+iu,\chi\right)\right|^{4}
\ll_{\epsilon}(Q(|u|+1))^{1+\epsilon}
\end{align*}
for $\epsilon>0$, where $S(Q)$ is the set of all real primitive characters
of conductor at most $Q$.
\end{lem}

\begin{lem}\label{lemCor} (see, e.g. Corollary 8.13, 8.19 of \cite{IwBook}) 
Let $a,b$ be real numbers with $b-a\geq1$ 
and $f(x)$ a real function. 
Suppose that $f(x)$ satisfies that 
$\Lambda_{2}\leq |f^{(2)}(x)| \leq \eta_{2}\Lambda_{2}$ and 
$\Lambda_{3}\leq |f^{(3)}(x)| \leq \eta_{3}\Lambda_{3}$ on $[a,b]$ 
for some $\Lambda_{2},\Lambda_{3}>0$ and $\eta_{2},\eta_{3}\geq1$. 
Then  we have 
\begin{align*}
\sum_{a<n<b}e(f(n))&\ll \eta_{2}\Lambda_{2}^{1/2}(b-a)+\Lambda_{2}^{-1/2},\\
\sum_{a<n<b}e(f(n))&\ll \eta_{3}^{1/2}\Lambda_{3}^{1/6}(b-a)+\Lambda_{3}^{-1/6}(b-a)^{1/2},
\end{align*}
where the implied constants are absolute.  
\end{lem}

\noi{\it Proof of Theorem \ref{thm}.} 
We first study $\psi_{s}(x)$. 
Denote by $s=\sigma+iT$ with $1/2<\sigma<1$ and 
define the value $\lambda_{q}(D)$ by $L(z,D)=\sum_{q\geq1}\lambda_{q}(D)q^{-z}$. 
Due to \S 5 of \cite{SY}, we have
\begin{align*} 
L(1,D)=\sum_{q\geq1}\lambda_{q}(D)q^{-1}e^{-q/U}
-\frac{1}{2\pi i}\int_{\frac{1}{2}-i\infty}^{\frac{1}{2}+i\infty}L(z,D)U^{z-1}\Gam(z-1)dz
\end{align*}
for $U>0$. 
Then $\psi_{s}(x)$ is written as follows.
\begin{align*}
\psi_{s}(x)=&2\sum_{q\geq1}q^{-1}e^{-q/U}
\sum_{3\leq t<X}\lambda_{q}(t^{2}-4)\left(1-\frac{\epsilon(t)^{2}}{x}\right)\epsilon(t)^{1-2s}\\
&-\frac{1}{\pi i}\int_{\frac{1}{2}-i\infty}^{\frac{1}{2}+i\infty}U^{z-1}\Gam(z-1)
\sum_{3\leq t<X}L(z,t^{2}-4)\left(1-\frac{\epsilon(t)^{2}}{x}\right)\epsilon(t)^{1-2s}dz\\
=:&A_{1}+A_{2}.
\end{align*}

By the definition of $L(z,D)$, we see that 
\begin{align*}
\sum_{3\leq t<X}\left|L\left(\frac{1}{2}+iu,t^{2}-4\right)\right|^{4}
\ll_{\epsilon} X^{\epsilon}\sum_{\chi\in S(X^{2})}\left| L\left(\frac{1}{2}+iu,\chi\right)\right|^{4},
\end{align*}
and then we get  
\begin{align*}
\sum_{3\leq t<X}\left|L\left(\frac{1}{2}+iu,t^{2}-4\right)\right| 
\ll_{\epsilon}  X^{\frac{5}{4}+\epsilon}(|u|+1)^{\frac{1}{4}+\epsilon}
\end{align*}
due to Lemma \ref{lemHB} and H\"{o}lder's inequality.
Thus $A_{2}$ can be bounded by 
\begin{align}\label{A2}
A_{2}&\ll_{\epsilon}U^{-\frac{1}{2}}\int_{\frac{1}{2}-i\infty}^{\frac{1}{2}+i\infty}|\Gam(z-1)|
\sum_{3\leq t<X}|L(z,t^{2}-4)|\epsilon(t)^{1-2\sigma}dz
\ll_{\epsilon} U^{-\frac{1}{2}}x^{\frac{9}{8}-\sigma+\epsilon}.
\end{align}

According to Lemma 2.3 of \cite{SY}, 
we have 
\begin{align}\label{A0}
A_{1}=&\sum_{q\geq1}q^{-1}e^{-q/U}\sum_{q_{1}^{2}q_{2}=q}q_{2}^{-1}\sum_{k\bmod{q_{2}}}S(k^{2},1;q_{2})
\sum_{3\leq t<X}e\left( \frac{kt}{q_{2}}-\frac{T}{\pi}\log{\epsilon(t)} \right)
\left(1-\frac{\epsilon(t)^{2}}{x} \right)\epsilon(t)^{1-2\sigma},
\end{align}
where $e(x):=e^{2\pi i x}$ and 
$$S(n,m;c):=\sum_{\begin{subarray}{c}a,b\bmod{c} \\ ab\equiv 1\bmod{c}\end{subarray}}
e\left(\frac{am+bn}{c}\right)$$ is the Kloostermann sum. 
It is easy to see that Lemma \ref{lemCor} 
yields the exponential sum estimate 
\begin{align*}
\sum_{N\leq t<2N}e\left( \frac{kt}{q_{2}}-\frac{T}{\pi}\log{\epsilon(t)} \right)
\ll \min(T^{\frac{1}{2}}+T^{-\frac{1}{2}}N,T^{\frac{1}{6}}N^{\frac{1}{2}})
\ll \begin{cases} T^{\frac{1}{6}}N^{\frac{1}{2}}, &(N<T^{2/3}),\\
T^{\frac{1}{2}}, & (T^{2/3}\leq N<T),\\ 
T^{-\frac{1}{2}}N, & (N\geq T)
\end{cases}
\end{align*}
for $N>1$, and then the sum over $t$ in \eqref{A0} can be bounded by 
\begin{align*}
\sum_{3\leq t<X}e\left( \frac{kt}{q_{2}}-\frac{T}{\pi}\log{\epsilon(t)} \right)
\left(1-\frac{\epsilon(t)^{2}}{x} \right)\epsilon(t)^{1-2\sigma}
\ll T^{-\frac{1}{2}}x^{1-\sigma}+T^{\frac{7}{6}-\frac{4}{3}\sigma}+T^{\frac{1}{6}}.
\end{align*}
Since $S(k^{2},1;q_{2})\ll_{\epsilon}q_{2}^{1/2+\epsilon}$, 
we can estimate $A_{1}$ as follows. 
\begin{align}\label{A1}
A_{1}\ll U^{1/2+\epsilon}(T^{-\frac{1}{2}}x^{1-\sigma}+T^{\frac{7}{6}-\frac{4}{3}\sigma}+T^{\frac{1}{6}}).
\end{align}

Combining \eqref{A2}, \eqref{A1} and Proposition \ref{prop1}, 
we have  
\begin{align*}
\frac{Z'_{\Gam}(s)}{Z_{\Gam}(s)}\ll_{\epsilon} 
T^{-\frac{1}{4}}x^{\frac{17}{16}-\sigma+\epsilon}
+(T^{\frac{7}{12}-\frac{2}{3}\sigma}+T^{\frac{1}{12}})
+x^{\frac{9}{16}-\frac{1}{2}\sigma+\epsilon}+T^{-2}x^{1-\sigma}+T^{1+\epsilon}x^{1/2-\sigma+\epsilon}.
\end{align*}
The result of Theorem \ref{thm} for $\frac{1}{2}<\sigma\leq \frac{5}{8}$ 
is derived from the above with $x=T^{\frac{20}{9}}$. 
We can obtain the result for $\frac{5}{8}<\sigma <1$ 
by the Phragmen-Lindel\"{o}f convexity theorem. \qed

\vcpt

\noi{\bf Acknowledgment.} 
The author was supported by JST CREST no.JPMJCR14D6 
and JSPS Grant-in-Aid for Scientific Research (C) 
no. 17K05181.

\end{document}